\documentclass[12pt,a4paper]{article}
\usepackage[utf8]{inputenc}
\usepackage[english]{babel}
\usepackage{amsmath,amsfonts,amssymb,amsthm}
\title{Modular forms on $\SU(2,1)$ with weight $\frac{1}{3}$}
\author{Eberhard Freitag and Richard M. Hill}
\date{}

\usepackage[english]{babel}
\usepackage[utf8]{inputenc}
\usepackage{amsmath,amssymb,amsfonts,amsthm}
\usepackage[dvipsnames]{xcolor}
\usepackage[colorlinks,linkcolor=Purple]{hyperref}

\usepackage[top=2.5cm, left=2cm, right=2cm, bottom=4.0cm]{geometry}

\newtheorem{definition}{Definition}
\newtheorem{corollary}{Corollary}

\newtheorem{lemma}{Lemma}
\newtheorem{proposition}{Proposition}
\newtheorem{theorem}{Theorem}

\newcommand{\optional}[1]{}

\DeclareMathOperator{\ab}{ab}

\DeclareMathOperator{\cts}{{cts}}

\DeclareMathOperator{\denom}{Denom}

\DeclareMathOperator{\GL}{GL}

\DeclareMathOperator{\Hom}{Hom}

\DeclareMathOperator{\meas}{meas}

\DeclareMathOperator{\Norm}{N}

\DeclareMathOperator{\PU}{PU}

\DeclareMathOperator{\SL}{SL}

\DeclareMathOperator{\SU}{SU}

\DeclareMathOperator{\Tr}{Tr}
\DeclareMathOperator{\U}{U}

\DeclareMathOperator{\verl}{Verl}

\newcommand{\C}{\mathbb{C}}
\newcommand{\F}{\mathbb{F}}

\newcommand{\Q}{\mathbb{Q}}

\newcommand{\Z}{\mathbb{Z}}

\newcommand{\cH}{\mathcal{H}}

\newcommand{\utau}{\underline{\tau}}

\begin{document}
\maketitle

\begin{abstract}
In this note, we describe several new examples of holomorphic modular forms on the group $\SU(2,1)$.
These forms are distinguished by having weight $\frac{1}{3}$.
We also describe a method for determining the levels at which one should expect to find such
fractional weight forms.
\end{abstract}

\tableofcontents

\section{Introduction}
\label{sec:intro}

In this note, we describe some new examples of holomorphic modular forms on the real Lie group $\SU(2,1)$.
The forms have weight $\frac{1}{3}$, and their levels are certain congruence subgroups of $\SU(2,1)$.
We recall that various half-integral weight forms on $\SU(2,1)$ have been studied in the past,
but our examples seem to be the first whose weight is not half-integral.
Indeed, they are not metaplectic forms, in the sense that they do not arise from
automorphic forms on the metaplectic cover of an ad\`ele group.
This is because the ad\`elic metaplectic covers of forms of $\SU(2,1)$ are at most two-fold covers.

\paragraph{Notation.}
To state our results a little more precisely, we introduce some notation.

Consider the Hermitian form of $\C^3$ defined by $\langle v,w\rangle = \bar v^t J w$, where $J=\begin{pmatrix} &&1\\&1\\1\end{pmatrix}$.
Here and throughout the paper, we write $\bar v^t$ for the conjugate transpose of $v$.
The group $\SU(2,1)$ consists of the matrices in $\SL_3(\C)$ which preserve this Hermitian form.
Equivalently, the elements of $	\SU(2,1)$ are the matrices $g$ in $\SL_3(\C)$ which satisfy $\bar g^t J g = J$.\footnote{One could of course replace $J$ by any
Hermitian matrix of signature $(2,1)$.
Our choice of $J$ is convenient, because the upper triangular matrices of $\SU(2,1)$ form a Borel subgroup  and the diagonal matrices form a maximal torus with this choice of $J$.
}

The symmetric space associated to this Lie group may be identified
with the following complex manifold:
\[
	\cH = \left\{ \utau = \begin{pmatrix}
		\tau_1\\ \tau_2
\end{pmatrix}	 \in \C^2 :
	\left\langle
	\begin{pmatrix} \tau_1 \\ \tau_2 \\ 1\end{pmatrix},\begin{pmatrix} \tau_1 \\ \tau_2  \\ 1\end{pmatrix}\right\rangle < 0 \right\}.
\]
Let $g=\begin{pmatrix} A&B\\C&D \end{pmatrix}$
be an element of $\SU(2,1)$,
where $A$ is a $2\times 2$ matrix, $D$ is a complex number, etc.
The action of $g$ on $\cH$ is given by
\[
	g*\utau
	=\frac{1}{C\utau+D} \cdot (A\utau+B).
\]
We also define
\begin{equation}
	\label{eqn: j definition}
	j(g,\utau)
	= C \utau + D.
\end{equation}
The function $j(g,\utau)$ satisfies the usual condition of a multiplier system:
\begin{equation}
	\label{eqn: j multiplier system}
	j(gh,\utau) = j(g,h * \utau)\cdot j(h,\utau),
	\qquad
	g,h \in \SU(2,1),\quad \utau \in \cH.
\end{equation}

\begin{definition}
	\label{defn:frac wt multiplier system}
	Let $\Gamma$ be an arithmetic subgroup of $\SU(2,1)$
	and let $\frac{a}{b}$ be a rational number.
	A function $\ell : \Gamma \times \cH \to \C^\times$ is called
	a \emph{fractional weight multiplier system} of weight $\frac{a}{b}$
	if it satisfies the following conditions:
	\begin{enumerate}
		\item
		$\ell$ is a multiplier system, i.e.
		$\ell(gh,\utau)= \ell(g,h*\utau)\cdot \ell(h,\utau)$
		for all $g,h\in \Gamma$ and $\utau \in \cH$;
		\item
		There exists a function
		$\chi : \Gamma \to \C^\times$ such that
		\(
			\ell(g,\utau)^b = \chi(g)\cdot j(g,\utau)^a
		\)
		for all $g\in \Gamma$ and $\utau \in \cH$.
		Such a function $\chi$
		must be a character of $\Gamma$
		by condition 1 and (\ref{eqn: j multiplier system}).

		\item
		For each $g\in \Gamma$ the function $\utau\mapsto \ell(g,\utau)$ is continuous (and hence holomorphic) on $\cH$.
	\end{enumerate}
\end{definition}

Half-integral weight multiplier systems have been known about for some time, and have been studied (for example) in \cite{wang williams} and \cite{jalal}.
The half-integral weight multiplier systems
arise from the theory automorphic forms on the metaplectic double cover of a unitary ad\`ele group.
However, no other fractional weight multiplier systems on $\SU(2,1)$ have been shown to exist until recently.

By a \emph{fractional weight modular form} of level $\Gamma$
with weight $\frac{a}{b}$, we shall mean a holomorphic function
$f: \cH \to \C$, such that for all $\utau\in \cH$ and all $\gamma\in \Gamma$, we have
\[
	f(\gamma*\utau) = \ell(\gamma,\utau) \cdot f(\utau),
\]
where $\ell$ is a multiplier system with weight $\frac{a}{b}$.
The usual growth condition in cusps is redundant for holomorphic modular forms on $\SU(2,1)$.

\paragraph{Statement of results.}

The purpose of this paper is to describe certain
modular forms on $\SU(2,1)$ of weight $\frac{1}{3}$.
The levels of these forms are congruence subgroups of
the following group
\begin{equation}
	\label{eqn:Gamma(1) defn}
	\Gamma(1) = \SU(2,1) \cap \SL_3(\Z[\zeta]),\qquad
	\text{where } \zeta = e^{\frac{2 \pi i}{3}}.
\end{equation}
(The reason for focusing on subgroups of $\Gamma(1)$
is because we make use a presentation for $\Gamma(1)$
found in \cite{falbel parker}
and the description of integral weight forms in \cite{freitag salvati}).
For a non-zero $\beta\in \Z[\zeta]$, we shall write $\Gamma(\beta)$
for the principal congruence subgroup of level $\beta$ in $\Gamma(1)$.
The levels of our forms are intermediate
groups between $\Gamma(\sqrt{-3})$ and $\Gamma(3)$, all of which have index $3$ in $\Gamma(\sqrt{-3})$.

One way in which one might expect to construct fractional weight forms is to use the theory of Borcherds products associated to the lattice $\Z[\zeta]^3$.
Recall that such a Borcherds product
is a meromorphic modular form $\Psi:\cH \to \C\cup\{\infty\}$, whose divisor
is a certain integer linear
combination of Heegner divisors.
Given a positive integer $n$ and a congruency class in $v \in \frac{1}{\sqrt{-3}}\Z[\zeta]^3 / \Z[\zeta]^3$, the Heegner divisor $D_{n,v}$ is defined by
\begin{equation}
	\label{eqn:Heegner}
	D_{n,v}
	=
	\begin{cases}
	\displaystyle{\sum_{w \in \Z[\zeta]^3 + v \;:\; \langle w, w \rangle = n}
	w^\perp} & \text{if $2v \in \Z[\zeta]^3$}\\[8mm]
	\displaystyle{\frac{1}{2}\sum_{w \in \Z[\zeta]^3 + v \;:\; \langle w, w \rangle = n}
	w^\perp} & \text{if $2v \not\in \Z[\zeta]^3$,}
	\end{cases}
\end{equation}
where $w^\perp$ is the set of $\utau \in \cH$ such that
$\langle \begin{pmatrix}\utau\\ 1\end{pmatrix},w\rangle = 0$.

Notice that if $w$ is a vector occurring in the sum (\ref{eqn:Heegner}), then $\zeta w$ and $\zeta^2 w$ also arise in the sum, and we have $w^\perp = (\zeta w)^\perp = (\zeta^2 w)^\perp$.
Hence each Heegner divisor $D_{n,v}$ is a multiple of $3$.
This implies that each Borcherds lift to $\SU(2,1)$ must be the cube of a meromorphic function on $\cH$.
If a Borcherds lift has weight $k$, then its cube root will be a fractional weight form of weight $\frac{k}{3}$, implying that there is a multiplier system with weight $\frac{k}{3}$.
The level of a Borcherds product associated with the lattice $\Z[\zeta]^3$ is the principal congruence subgroup $\Gamma(\sqrt{-3})$.
In view of this, it makes sense to begin searching for
third-integral weight forms at level $\Gamma(\sqrt{-3})$.
In spite of this optimism, we obtain the following result:

\begin{theorem}
	Every multiplier system on the group $\Gamma(\sqrt{-3})$
	has integral weight.
	Hence every modular form of level $\Gamma(\sqrt{-3})$ has integral weight.
\end{theorem}

This result is proved in \autoref{thm:wt denominator}; it is a proof by contradiction, using a carefully chosen relation in the group $\Gamma(\sqrt{-3})$.
Note that the theorem implies that the weight of each of the Borcherds product associated to the lattice $\Z[\zeta]^3$ must be a multiple of $3$.
Searching a little further, we do indeed find some fractional weight multiplier systems:

\begin{theorem}
	\label{thm:intro 2}
	There are thirteen subgroups $\Gamma$ of index $3$ in $\Gamma(\sqrt{-3})$
	and containing $\Gamma(3)$, such that $\Gamma$ has a multiplier system of weight $\frac{1}{3}$.
\end{theorem}

The result is proved computationally, by finding a presentation for each of the subgroups.
The thirteen subgroups mentioned in the theorem are listed in
\autoref{thm:list}.

\autoref{thm:intro 2} gives us a list of levels where we might conceivably
find a modular form of weight $\frac{1}{3}$.
Using results of \cite{freitag salvati}, we obtain a description of the ring of
integral weight modular forms of level $\Gamma(3)$, as well as the action of $\Gamma(\sqrt{-3})/\Gamma(3)$ on this ring.
We find a list of forms of level $\Gamma(3)$ and weight $1$, all of whose divisors are multiples of $3$.
By considering the action of $\Gamma(\sqrt{-3})/\Gamma(3)$ on these forms, we discover that
each of them is a modular form (with character) of some larger level.
Each of these larger levels is one of the groups found in Theorem 2.
The cube roots of these forms are forms of weight $\frac{1}{3}$.
In this way, we prove the following:

\begin{theorem}
	\label{thm:intro 3}
	For twelve of the thirteen groups $\Gamma$ listed in \autoref{thm:list}, there exists a holomorphic modular form of weight $\frac{1}{3}$ and level $\Gamma$.
\end{theorem}

\paragraph{Further comments and questions.}

\autoref{thm:intro 3}
 leaves an obvious open question: why is there no form of weight $\frac{1}{3}$ on the remaining group $\Gamma$ in our list.
We have no clear answer to this question.

We have investigated the multiplier systems
on many other congruence subgroups of $\Gamma(1)$.
For each congruence subgroup that we have investigated,
the denominator of the weight of every multiplier system is a factor of $6$.
It is tempting to think that this might be true
for all congruence subgroups of $\Gamma(1)$.
One might even guess that if we replace the field $\Q(\zeta)$ involved in the construction of $\Gamma(1)$ by another CM field $k$, then a similar bound might be the number of roots of unity in $k$.
However, we have not proved any such result;
neither do we have any numerical evidence beyond the case $k=\Q(\zeta)$.

The situation for non-congruence subgroups is very different.
For every positive integer $n$,
there exists a (non-congruence) subgroup $\Gamma_n$ of finite index in $\Gamma(1)$,
such that $\Gamma_n$ has a multiplier system of weight $\frac{1}{n}$.
This fact follows from a more general result
proved in \cite{hill res-fin},
and independently in \cite{Stover Toledo}.

For $d > 2$, it is not known whether there exist
any modular forms whose level is a congruence subgroup $\SU(d,1)$, and whose weight has denominator greater than $2$.
However, such forms certainly do exist at non-congruence levels by results in \cite{hill res-fin}.

\paragraph{Organization of the paper.}
The paper is organized as follows.
Suppose $\Gamma$ is an arithmetic subgroup of $\SU(2,1)$.
We prove (see \autoref{thm:finiteness}) that there is an upper bound on the denominators of weights of multiplier systems on $\Gamma$.
We define in \autoref{section:background}
a natural number called the
\emph{weight denominator} of $\Gamma$.
The weight denominator
is the lowest common denominator of the weights of all multiplier systems on $\Gamma$.

Let $\tilde\Gamma$ be the pre-image of $\Gamma$ in the universal cover of $\SU(2,1)$.
One may easily calculate the weight denominator of $\Gamma$ if one has a presentation of $\tilde\Gamma$.
In \autoref{sec:presentations}, we recall some methods for obtaining such presentations. These methods give us a computational approach to determining the weight denominator of a given arithmetic group.

In \autoref{sec:results} we describe in detail how the
methods from \autoref{sec:presentations} have been implemented for
subgroups of finite index in the group $\Gamma(1)$.
We calculate the weight denominators of various arithmetic groups, and we list in \autoref{sec:index 3 subgroups} thirteen congruence subgroups with weight denominator $3$, all with index $3$ in $\Gamma(\sqrt{-3})$.

In \autoref{section: modular forms} we construct,
for twelve of the thirteen groups found in \autoref{sec:index 3 subgroups},
a modular form of weight $\frac{1}{3}$.

\section{Background on SU(2,1)}
\label{section:background}

In this section we shall define and investigate a positive integer, which we call the \emph{weight denominator} of an arithmetic subgroup of $\SU(2,1)$.

Let $\widetilde{\SU(2,1)}$ be the universal cover of $\SU(2,1)$.
For an arithmetic group $\Gamma$ we write $\tilde\Gamma$ for the
pre-image of $\Gamma$ in $\widetilde{\SU(2,1)}$.
The \emph{weight denominator}
of $\Gamma$ is defined to be the
order of the kernel of the
projection map
$\tilde\Gamma / [\tilde\Gamma,\tilde\Gamma]
\to \Gamma / [\Gamma,\Gamma]$.
We prove in \autoref{thm:finiteness}
that the weight denominator is finite.
We prove in \autoref{cor:nmax} that
 $\Gamma$ has a multiplier system of rational weight $w$ if and only if
the denominator of $w$ is a factor of the weight denominator of $\Gamma$.
Some other simple properties of the weight denominator are proved in \autoref{prop:wt denominator facts}.
Methods for computing the weight denominator
of an arithmetic group are discussed
in \autoref{sec:presentations} and
\autoref{sec:results}.

\subsection{The universal cover}

The Lie group $\SU(2,1)$ has fundamental group $\Z$.
We shall write $\widetilde{\SU(2,1)}$ for its universal cover.
We therefore have a central extension of groups
\[
	1 \to \Z \to \widetilde{\SU(2,1)} \to \SU(2,1) \to 1.
\]
Such extensions are classified by
elements of the measurable cohomology group
$H^2_{\meas}(\SU(2,1),\Z)$ (see \cite{Moore}).
In order to perform calculations, it will be helpful to have a specific measurable 2-cocycle $\sigma$
representing this group extension.
We describe such a cocycle in \autoref{prop:sigma} below. As a first step, we define a function $X : \SU(2,1)\to \C$ by
\[
	X\begin{pmatrix}
		*&*&*\\
		*&*&*\\
		a&b&c
	\end{pmatrix}
	= \begin{cases}
		-a & \text{if $a\ne 0$,}\\
		c & \text{if $a = 0$.}
	\end{cases}
\]

\begin{lemma}
	\label{lemma:half-plane}
	Let $g\in \SU(2,1)$ and $\utau \in \cH$.
	Then $X(g) \ne 0$ and
	the complex number $\frac{j(g,\utau)}{X(g)}$ has positive real part, where $j(g,\utau)$ is defined in (\ref{eqn: j definition}).
\end{lemma}

\begin{proof}
	Let $\begin{pmatrix}a&b&c\end{pmatrix}$ be the bottom row of the matrix $g$.
	Since $g\in \SU(2,1)$ we must have $\bar a c + \bar b b + \bar c a =0$.
	We shall divide the proof into two cases depending on whether or not $a=0$.
	
	Consider first the case where $a=0$.
	The equation above implies $b=0$, and therefore $c \ne 0$.
	In such cases $X(g)=j(g,\utau)=c$ for all $\utau\in \cH$ so the lemma is true in this case.
	
	Suppose from now on that $a\ne 0$.
	In this case we have $X(g)=-a$, and in particular $X(g) \ne 0$.
	The group $\SU(2,1)$ contains the following torus:
	\[
		T = \{ t_z : z\in \C^\times \},
		\qquad
		\text{where }
		t_z = \begin{pmatrix}
			\frac{1}{\bar z} \\ & \frac{\bar z}{z} \\&& z
		\end{pmatrix}.
	\]
	It is trivial to check that $X(t_z\cdot  g) = z \cdot X(g)$
	and $j(t_z \cdot g , \utau) = z\cdot j(g,\utau)$.
	In view of this, it is sufficient to prove the lemma in the case $a=1$, so that we have $2\Re(c) = - |b|^2$.
	We must check in this case that $j(g,\utau)$ has negative real part.
	
	Let $\utau = \begin{pmatrix} \tau_1 \\ \tau_2\end{pmatrix}$.
	Since $\utau\in \cH$ we must have
	$2\Re(\tau_1) + |\tau_2|^2 < 0$.
	This implies:
	\[
		\Re(j(g,\utau))
		\;=\; \Re( \tau_1  + b \tau_2 + c)
		\;< \;-\frac{|\tau_2|^2}{2} + \Re (b \tau_2)  - \frac{|b|^2}{2}
		\;=\; -\left|\tau_2-\bar b\right|^2
		\; \le\; 0.
	\]
\end{proof}

\autoref{lemma:half-plane} allows us to define, for each $g\in \SU(2,1)$, a branch $\tilde j(g,-)$ of the logarithm of $j(g,-)$ as follows:
\begin{equation}
	\label{eqn:j-tilde}
	\tilde j(g,\utau) = \log\left( \frac{j(g,\utau)}{X(g)} \right) + \log(X(g)),
\end{equation}
where each logarithm in the right hand side
of (\ref{eqn:j-tilde}) is defined to be continuous away from the negative real axis and satisfy $-	\pi < \Im ( \log z) \le \pi$.
For each fixed $g\in\SU(2,1)$ the function $\tilde j(g,\utau)$ is continuous in $\utau$ by \autoref{lemma:half-plane}.

The multiplier system condition (\ref{eqn: j multiplier system}) on $j(g,\utau)$ implies the following congruence for $\tilde j$:
\[
	\tilde j(gh,\utau) \equiv \tilde j(g,h*\utau) + \tilde j(h,\utau) \mod 2\pi i\cdot \Z.
\]
In particular, we may define for $g,h \in\SU(2,1)$ an integer $\sigma(g,h)$ by
\begin{equation}
	\label{eqn:sigma definition}
	\sigma(g,h) = \frac{1}{2\pi i} \Big(\tilde j(gh,\utau) - \tilde j(g,h*\utau) - \tilde j(h,\utau)\Big).
\end{equation}
The right hand side of (\ref{eqn:sigma definition}) is independent of $\utau$, since it is a continuous $\Z$-valued function of $\utau$.

\begin{proposition}
	\label{prop:sigma}
	The function $\sigma$ defined in (\ref{eqn:sigma definition}) is an inhomogeneous measurable 2-cocycle, whose cohomology class in $H^2_{\meas}(\SU(2,1),\Z)$
	corresponds to the universal cover of $\SU(2,1)$.
\end{proposition}

\begin{proof}
	The function $\sigma$ is evidently measurable, and it is a short exercise
	using (\ref{eqn:sigma definition}) to verify the 2-cocycle relation:
	\[
		\sigma(g_1,g_2)+\sigma(g_1g_2,g_3)
		=\sigma(g_1,g_2g_3)+\sigma(g_2,g_3).
	\]
	By Calvin Moore's theory of measurable cohomology
	\cite{Moore},
	there is a central extension of Lie groups
	corresponding to $\sigma$.
	It remains to show that this extension is the universal cover.
	For the moment, we'll write $\widetilde{\SU(2,1)}$ for the
	extension of $\SU(2,1)$ corresponding to $\sigma$.
	Explicitly, $\widetilde{\SU(2,1)}$ is the set $\SU(2,1) \times \Z$, with
	the group operation given by
	\begin{equation}
		\label{eqn:sigma multiplication}
		(g,n)\cdot (g',n')
		=(gg', n+n'+\sigma(g,g')).
	\end{equation}
	To prove that $\widetilde{\SU(2,1)}$ is the universal cover of $\SU(2,1)$, it's sufficient to prove that $\widetilde{\SU(2,1)}$ is connected.
	(Note that the topology on the Lie group $\widetilde{\SU(2,1)}$ is not the
	product topology on $\SU(2,1) \times \Z$;
	however the Borel measurable subsets of the Lie group $\widetilde{\SU(2,1)}$
	coincide with the Borel measurable subsets of the product $\SU(2,1)\times \Z$).
	
	Let $\tilde T$ be the pre-image in $\widetilde{\SU(2,1)}$ of the following torus in $\SU(2,1)$:
	\[
		T = \left\{t_z:
		z\in \C^\times\right\},
		\qquad \text{where }
		t_z=\begin{pmatrix} \frac{1}{\bar z} \\& \frac{\bar z}{z}\\ && z \end{pmatrix}.
	\]
	Note that on $T$, we have simply $\tilde j(t_z,\utau) = \log(z)$, where as before we are taking $-\pi < \Im (\log z) \le \pi$.
	This easily implies that there is an
	isomorphism $\Psi:\tilde T \to \C$ given by
	\begin{equation}
		\label{eqn:T tilde isomorphism}
		\Psi(t_z,n) = \log(z) - 2\pi i\cdot n.
	\end{equation}
	The map $\Psi$ and its inverse are measurable.
	Since measurable homomorphisms of Lie groups are continuous, it follows that $\Psi$ is a homeomorphism.
	In particular $\tilde T$ is connected.
	Hence all elements of the subgroup $\Z \subset \widetilde{\SU(2,1)}$ are connected by paths in $\tilde T$ to the identity element.
	This implies that $\widetilde{\SU(2,1)}$ is connected.
\end{proof}

From now on, we shall refer to elements of $\widetilde{\SU(2,1)}$ as pairs $(g,n) \in \SU(2,1)\times \Z$, with the group operation given by
(\ref{eqn:sigma multiplication}).
To find out whether a fractional weight multiplier system exists on a group $\Gamma$, we shall use the following result.

\begin{proposition}
	\label{prop:multiplier character equivalence}
	Let $\Gamma$ be an arithmetic subgroup of $\SU(2,1)$; let $\tilde\Gamma$ be the pre-image of $\Gamma$ in $\widetilde{\SU(2,1)}$, and let $w$ be a rational number.
	
	For any multiplier system $\ell$ on $\Gamma$  of weight $w$,
	there is a character $\Phi: \tilde \Gamma \to \C^\times$ defined by
	\begin{equation}
		\label{eqn:Phi defn}
		\Phi(g,n)
		=
		\ell(g,\utau) \cdot\exp\Big(w\left( 2 \pi i \cdot n - \tilde j(g,\utau)\right)\Big).
	\end{equation}
	In particular, the right hand side of (\ref{eqn:Phi defn}) is independent of $\utau\in \cH$.
	
	Conversely, given any character $\Phi:\tilde\Gamma \to \C^\times$, equation (\ref{eqn:Phi defn})
	defines a weight $w$ multiplier system $\ell(g,\utau)$ for any
	rational number $w$ satisfying
	$\Phi(I_3,1)=e^{2\pi i\cdot w}$.
\end{proposition}

\begin{proof}
	Let $\ell(g,\utau)$ be a multiplier system on $\Gamma$ with weight $w=\frac{a}{b}$ and character $\chi$.
	If we fix $(g,n)\in\tilde\Gamma$, then
	our formula (\ref{eqn:Phi defn}) for $\Phi(g,n)$ is a continuous function of
	$\utau\in \cH$.
	Furthermore, $\Phi(g,n)^b = \chi(g)$, which does not depend on $\utau$.
	Hence $\Phi(g,n)$ does not depend on $\utau$.
	To show that $\Phi$ is a homomorphism,
	we calculate as follows:
	\begin{align*}
		\Phi\big((g,n)\cdot (h,m)\big)
		&=\Phi\big(gh,n+m +\sigma(g,h)\big)\\
		&=\ell(gh,\utau)\cdot
		\exp\Big(w\big( 2 \pi i (n+m +\sigma(g,h)) - \tilde j(gh,\utau)\big)\Big)
	\end{align*}
	Substituting the definition (\ref{eqn:sigma definition}) of $\sigma$, we get
	\begin{align*}
		\Phi\big((g,n)\cdot (h,m)\big)
		&=\ell(gh,\utau)\cdot
		\exp\Big(w\big( 2 \pi i (n+m)
		- \tilde j(g,h*\utau)
		- \tilde j(h,\utau)
		\big)\Big).
	\end{align*}
	Using the multiplier system property of $\ell(g,\utau)$, we have
	\begin{align*}
		\Phi\big((g,n)\cdot (h,m)\big)
		&=\ell(g,h*\utau)\cdot \ell(h,\utau)\cdot
		\exp\Big(w\big( 2 \pi i (n+m) - \tilde j(g,h*\utau)-\tilde j(h,\utau) \big)\Big)\\
		&=\Phi(g,n) \cdot \Phi(h,m).
	\end{align*}
	
	Conversely, suppose $\Phi : \tilde \Gamma \to \C^\times$ is a character and $\Phi(I_3,1)=e^{2 \pi i \cdot w}$.
	The argument above may be reversed, to show that the function $\ell(g,\utau)=\Phi(g,0)\exp(w\cdot \tilde j(g,\utau))$ satisfies the multiplier system relation.
	Assuming $w=\frac{a}{b}$, we have
	\[
		\ell(g,\utau)^b
		=\Phi(g,0)^b
		\cdot j(g,\utau)^{a}.
	\]
	Therefore $\ell$ is a weight $w$ multiplier system with character $\chi(g) = \Phi(g,0)^b$.
\end{proof}

\subsection{The weight denominator}
\label{sec:wt denominator}
\begin{definition}
	Let $\Gamma$ be an arithmetic subgroup of $\SU(2,1)$ and let $\tilde{\Gamma}$ be the preimage of $\Gamma$ in $\widetilde{\SU(2,1)}$.
	Furthermore, let $\tilde\Gamma^{\ab}=\tilde{\Gamma}/[\tilde{\Gamma},\tilde{\Gamma}]$ be the
	abelianization of $\tilde\Gamma$.
	We define the \emph{weight denominator} $\denom(\Gamma)$ to be the order of the element $(I_3,1)$ in $\tilde{\Gamma}^{\ab}$.
\end{definition}

The following proposition is useful to know, although we do not use it in the rest of this paper.

\begin{theorem}
	\label{thm:finiteness}
	For any arithmetic subgroup $\Gamma$ of $\SU(2,1)$, the weight denominator is finite.
\end{theorem}

In the proof we'll use the following facts about the cohomology of a connected Lie group $G$:
\begin{itemize}
\item
There is an isomorphism $H^2_{\meas}(G,\Z) \cong \Hom(\pi_1(G),\Z)$
(see \cite{Moore}).
In particular $H^2_{\meas}(G,\Z)$ is torsion-free, so in injects into $H^2_{\meas}(G,\Z)\otimes\C$.
\item
By \cite{Wigner}, there is
an isomorphism
$H^\bullet_{\meas}(G,\Z) \otimes \C \cong H^\bullet_{\cts}(G,\C)$,
where $H^\bullet_{\cts}$ denotes continuous cohomology.
\item
If $\Gamma$ is a cocompact arithmetic subgroup of $G$ then by results in \cite{BorelWallach},
the restriction maps
$H^\bullet_{\cts}(G,\C) \to H^\bullet(\Gamma,\C)$ are all in injective.
\end{itemize}
Combining these results, we see that if
$\Gamma$ is cocompact in $G$, then the map
$H^2_{\meas}(G,\Z) \to H^2(\Gamma,\C)$ is injective.

\begin{proof}
	Let us suppose, for the sake of argument, that
	$\denom(\Gamma)$ is infinite.
	Hence the intersection of $\Z$ with $[\tilde\Gamma,\tilde{\Gamma}]$ is trivial,
	so $\Z$ injects into the finitely generated abelian group $\tilde{\Gamma}^{\ab}$.
	Choose a torsion-free subgroup $L$ of finite index in $\tilde\Gamma^{\ab}$ containing $\Z$.
	The map $\Z \to L$ has a left inverse $\phi : L \to \Z$.
	We may inflate $\phi$ to a map $\tilde\Gamma' \to \Z$, where $\tilde\Gamma'$ is the pre-image of
	$L$ in $\tilde{\Gamma}$.
	This implies $\tilde{\Gamma}' \cong \Gamma' \oplus \Z$, where $\Gamma'$ is the image of $\tilde\Gamma'$ in $\Gamma$.
	Hence the restriction of $\sigma$ to $\Gamma'$ is a coboundary.
	We'll use this assertion to obtain a contradiction by proving that the image of $\sigma$
	in $H^2(\Gamma',\C)$ is non-zero.

	In the case that $\Gamma'$ is cocompact in $\SU(2,1)$, the discussion above shows that the
	image of $\sigma$ in $H^2(\Gamma',\C)$ is non-zero, and we are done.
	In the case that $\Gamma'$ has cusps we need to be a little more careful.
	In this case, we may choose a Lie subgroup $G\subset\SU(2,1)$ isomorphic to $\SU(1,1)$, such that $\Gamma' \cap G$ is
cocompact in $G$ (see page 590-591 of \cite{Shimura}).
	The subgroup $G$ contains a conjugate of
	the matrix
	\[
		t_{-1} = \begin{pmatrix}
			-1 \\ & 1 \\ && -1
		\end{pmatrix}.
	\]
	Although $t_{-1}$ has order $2$, every pre-image of $t_{-1}$
	has infinite order in $\widetilde{\SU(2,1)}$
	by the isomorphism $\tilde T \cong \C$	in
	(\ref{eqn:T tilde isomorphism}).
	This shows that $\sigma$ does not split on $G$.
	Again by the discussion above, the image of $\sigma$ in
	$H^2(\Gamma'\cap G ,\C)$ is non-zero.
	The map $H^2_{\meas}(\SU(2,1),\Z) \to H^2(G\cap \Gamma', \C)$ factors as follows:
	\[
		H^2_{\meas}(\SU(2,1),\Z) \to H^2(\Gamma',\C) \to H^2 (G \cap \Gamma',\C).
	\]
	Hence the image of $\sigma$ in $H^2(\Gamma',\C)$ is non-zero.
\end{proof}

\begin{corollary}
	\label{cor:nmax}
	Let $\Gamma$ be an arithmetic subgroup of $\SU(2,1)$ and let $w$ be a rational number.
	There exists a weight $w$ multiplier system on $\Gamma$
	if and only if the denominator of $w$ is a factor of $\denom(\Gamma)$.
\end{corollary}

\begin{proof}
	By \autoref{prop:multiplier character equivalence}, the existence of a weight $w$ multiplier system is equivalent to the existence
	of a character $\Phi:\tilde\Gamma \to \C^\times$ satisfying $\Phi(I_3,1)= e^{2 \pi i \cdot w}$.
	Any such character $\Phi$ would factor as
	a map $\Phi:\tilde\Gamma^{\ab} \to \C^\times$.
	By Pontryagin duality,
	there exists a homomorphism $\Phi : \tilde \Gamma^{\ab} \to \C^\times$
	satisfying $\Phi(I_3,1)=e^{2 \pi i \cdot w}$ if and only if the denominator of $w$
	is a factor of the order of $(I_3,1)$ in
	$\tilde\Gamma^{\ab}$.
\end{proof}

\begin{proposition}
	\label{prop:wt denominator facts}
	Let $\Gamma$ and $\Gamma'$ be arithmetic subgroups of $\SU(2,1)$ with $\Gamma' \subset\Gamma$,
	and let $Z$ be the centre of $\SU(2,1)$.
	\begin{enumerate}
		\item
		$\denom(\Gamma)$ is a factor of $\denom(\Gamma')$.
		\item
		$\denom(\Gamma')$ is a factor of $n \cdot \denom(\Gamma)$, where $n$ is the index of $\Gamma'$ in $\Gamma$.
		\item
		$\denom(\Gamma \cdot Z) = \denom(\Gamma)$.
	\end{enumerate}
\end{proposition}

\begin{proof}
\begin{enumerate}
	\item
	By definition, $\denom(\Gamma)$ is the smallest positive integer $n$, such that $(I_3,n) \in [\tilde{\Gamma},\tilde{\Gamma}]$.
	The first statement follows because $[\tilde\Gamma',\tilde\Gamma']$ is a subgroup of $[\tilde\Gamma,\tilde\Gamma]$.
	\item
	Recall that if $H$ is a subgroup of finite index in a group $G$, then there is a transfer
	(or Verlagerung) homomorphism $\verl : G^{\ab} \to H^{\ab}$.
	For elements $z$ in the centre of $G$, the transfer map is given by $\verl(z)=z^n$, where $n$ is the index of $H$ in $G$.
	
	We therefore have a homomorphism
	$\verl :\tilde \Gamma^{\ab} \to \tilde{\Gamma}^{'\ab}$,
	which takes $(I_3,1)$ to $(I_3,n)$.
	Since $(I_3,\denom(\Gamma))$ is the identity
	element in $\tilde{\Gamma}^{\ab}$, it follows that $(I_3,n\cdot \denom(\Gamma))$ is the identity element in $\tilde{\Gamma}^{'\ab}$.
	\item
	Let $\tilde Z$ be the pre-image of $Z$ in $\widetilde{\SU(2,1)}$.
	By a general result on covering groups of Lie groups, $\tilde Z$ is the centre of $\widetilde{\SU(2,1)}$.
	Hence the commutator of an element of $\tilde Z$ with any other group element is trivial.
	This implies $[\widetilde{\Gamma \cdot Z}, \widetilde{\Gamma \cdot Z}] = [\tilde \Gamma, \tilde \Gamma]$, from which our result
	immediately follows.
\end{enumerate}
\end{proof}

\section{Calculus of presentations}
\label{sec:presentations}

To calculate the weight denominators of arithmetic groups $\Gamma$, we need some methods for producing presentations of the groups $\tilde\Gamma$.
In this section, we shall describe some methods for doing this.
We have used \cite{DLJohnson} as our main source for this section.
The methods described here apply to
arbitrary finitely presented groups.

\subsection{Covering groups}
\label{sec:covering groups}
Suppose $\tilde G$ is a central extension of a group $G$ by a cyclic group $\langle z \rangle$
of the form
\[
	1 \to \langle z \rangle \to \tilde G \to G \to 1.
\]
Assume that we have a presentation of $G$:
\[
	G = \langle g_1,\ldots,g_n | r_1 , \ldots , r_m \rangle,
\]
where each relation $r_i$ is a word in the generators $g_i$.
We shall describe a presentation of $\tilde G$.

For each generator $g_i$ in $G$, we choose a pre-image $\hat g_i$ of $g_i$ in $\tilde G$.
The group $\tilde G$ is generated by the elements
\[
	\hat g_1 , \ldots , \hat g_n , z.
\]
There are three obvious kinds of relation in
 $\tilde G$:
\begin{enumerate}
\item
For each relation $r_i$ in our presentation of $G$, we let $\hat r_i$ be the word in the generators of $\tilde G$, obtained by replacing each $g_i$ with the corresponding element $\hat g_i$ in $\tilde G$.
Since $r_i=1$ in $G$, it follows that $\hat r_i \in \langle z\rangle$, so we have a relation
$\hat r_i = z^{c_i}$ in $\tilde G$.

\item
Since the extension is assumed to be central,
the element $z$ is in the centre of $\tilde G$,
so we have a relation $[z,\hat g_i]$ for each generator $g_i$ in $G$.

\item
If $z$ has finite order $n$, then we have the relation $z^n=1$.
\end{enumerate}
The generators and relations listed above give a presentation of $\tilde G$.

\subsection{Subgroups}
\label{sec:Reidemeister Schreier}
Let $H$ be a subgroup of finite index in $G$.
Assume that we have a finite presentation of $G$
with generators $g_i$.
The method for constructing a finite presentation of $H$ is called the Reidemeister -- Schreier algorithm.
This algorithm constructs a certain directed graph, whose vertex set $R$ is a set of coset representatives for $H$ in $G$.

The Reidemeister--Schreier graph is constructed as follows:
\begin{enumerate}
	\item
	We begin by setting $R_0=R=\{1_G\}$, where $1_G$ is the identity element in the group $G$. 	\item
	Choose an element $r \in R_0$ and remove it from the set $R_0$.
	\item
	For each element $g$ which is either one of the generators $g_i$ of $G$ or its inverse $g_i^{- 1}$,
	consider the coset $Hrg$.
	Check whether $Hrg$ is one of the cosets $Hr'$ for some element $r'$ in $R$.
	\begin{enumerate}
		\item
		If $Hrg=Hr'$,
		then we must have $rg=hr'$ for some $h\in H$.
		We then add an edge from $r$ to $r'$,
		and label that edge $(g, h)$.
		\item
		If $Hrg$ is not one of these cosets,
		then we add $rg$ to our set $R$ and also to our set $R_0$.
		We also add an edge from $r$ to $rg$
		in our graph, labelled $(g,1)$.
	\end{enumerate}
	\item
	If $R_0$ is non-empty, then go back to step 2.
	Otherwise stop.
\end{enumerate}
The algorithm will terminate with a graph
whose vertices are labelled by
representatives $r$ for cosets $Hr$.
For each element $g$ which is either a generator $g_i$
or its inverse $g_i^{-1}$, there is an edge
from $r$ to a vertex $r'$, where $Hr g = Hr'$.
This edge is labelled $(g,h)$, where
the element $h\in H$ satisfies $rg=hr'$.
We shall refer to $g$ and $h$ as the $G$-label and $H$-label of the edge respectively.

We'll now explain how one uses the Reidemeister--Schreier graph to construct a finite presentation of $H$ from one of $G$.
Suppose we have a word $w = w_1 \cdots w_s$
in the generators $g_i^{\pm 1}$ of $G$.
For any vertex $r$ in the graph, we may form a path
starting at $r$, and moving along the vertices
with $G$-labels $w_1, w_2, \ldots, w_s$ in turn.
The end of the path will be
a vertex $r'$ such that
$Hr w_1 \cdots w_s = Hr'$.
If $h_1, \cdots, h_s$ are the $H$-labels of the edges in the path then we have
\[
	r\cdot w_1 \cdots w_s = h_s \cdots h_1\cdot  r'.
\]

If $w_1 \cdots w_s=1$ is one of our relations in $G$,
then this path will be a closed loop, i.e. $r'=r$.
We will therefore have
\[
	r =  r \cdot w_1 \cdots w_s = h_s \cdots h_1 \cdot r,
\]
This gives us a relation $h_s \cdots h_1 = 1$ in the subgroup $H$.

It is known that the subgroup $H$ is generated by
the elements $h$ which are the $H$-labels of edges in the Reidemeister--Schreier graph.
Furthermore, the set of all relations arising as described above gives a presentation of the subgroup $H$.

\subsection{Rewriting Presentations}
\label{sec:rewriting}
Typically, the Reidemeister--Schreier algorithm produces a presentation with many generators, when far fewer generators are actually needed.
The \emph{rewriting} process allows us to reduce the number of generators.

Suppose that we have a group presentation:
\[
	G = \langle g_1,\ldots,g_r | r_1,\ldots,r_t \rangle,
\]
and we also have a second set of generators
$\{h_1,\ldots,h_s\}$ for $G$.
We'll explain now how to find a presentation using $h_1,\ldots, h_s$ as generators.
This process is called ``rewriting'' the presentation.

\begin{enumerate}
	\item
	Find an expression for each element $h_i$ as
	a word in the generators $g_j$:
	\[
		h_i = w_i(g_1,\ldots,g_r).
	\]
	Then we can the add the generators $h_i$ to the presentation of $G$,
	together with the relations $h_i=w_i(g_1,\ldots,g_r)$.

	\item
	Find an expression for each $g_i$
	as a word in the generators $h_j$:
	\[
		g_i = x_i(h_1,\ldots,h_s).
	\]
	We can add the relations $g_i = x_i(h_1,\ldots,h_s)$ to the presentation of $G$.
	
	\item
	In each relation apart from
	the relation  $g_i = x_i(h_1,\ldots,h_s)$,
	we replace each occurrence of $g_i$ by the word $x_i(h_1,\ldots,h_s)$.

	\item
	Finally, we remove the generators $g_i$ and the relations $g_i=x_i(h_1,\ldots,h_s)$.
\end{enumerate}

\section{Weight denominators of certain arithmetic groups}

\label{sec:results}

This section concerns computer calculations.
The code for these calculations runs on sage version 9.0, and is available at \cite{hill code}.

As before (see (\ref{eqn:Gamma(1) defn})), we shall write $\Gamma(1)$ for the subgroup of $\SU(2,1)$ consisting of matrices whose entries are in the ring $\Z[\zeta]$
of Eisenstein integers, where $\zeta=e^{2\pi i /3}$.
Let $\Gamma(\sqrt{-3})$ denote the principal congruence subgroup in $\Gamma(1)$ of level $\sqrt{-3}$.
As mentioned in the introduction,
the group $\Gamma(\sqrt{-3})$ is a level where one might expect to find forms of third-integral weight.
In spite of this, we show in
\autoref{thm:wt denominator} that $\Gamma(\sqrt{-3})$ has weight denominator $1$.

Looking just a little further up the tower of congruence subgroups, we do indeed find groups with weight denominator $3$.
The group $\Gamma(3)$ has weight denominator $3$.
Furthermore, there are precisely 13 subgroups of index $3$ in $\Gamma(\sqrt{-3})$ containing $\Gamma(3)$, which have weight denominator $3$
(see \autoref{sec:index 3 subgroups}).

\subsection{The subgroup $\Upsilon$}

Rather than dealing with the group $\Gamma(\sqrt{-3})$ directly,
it will be slightly more convenient to
consider a certain subgroup $\Upsilon$ of index $3$ in $\Gamma(\sqrt{-3})$.
This is because $\Upsilon$ has a slightly smaller presentation (see \autoref{prop:Upsilon presentation} below).


\begin{lemma}
	\label{lem:upsilon direct sum}
	The group $\Gamma(\sqrt{-3})$ decomposes as a direct sum:
	$\Gamma(\sqrt{-3}) = \Upsilon  \oplus Z$,
	where $Z$ is the centre of $\SU(2,1)$ generated by $\zeta \cdot I_3$,
	and $\Upsilon $ is a subgroup of index $3$ defined as follows:
	\[
		\Upsilon
		=
		\{ (g_{i,j}) \in \Gamma(\sqrt{-3}) : g_{1,1} \equiv 1 \bmod 3 \}.
	\]
	If $(g_{i,j})$ is any element of $\Upsilon$
	then we have $g_{1,1} \equiv g_{2,2} \equiv g_{3,3} \equiv 1 \bmod 3$.
\end{lemma}

Note that by \autoref{prop:wt denominator facts} the groups $\Upsilon$ and $\Gamma(\sqrt{-3})$ have the same weight denominator, since they are the same modulo the centre.

\begin{proof}
	Let $g= (g_{i,j})$ be an element of $\Gamma(\sqrt{-3})$.
	If $v$ and $w$ are the first and third columns of $g$, then we must have $\langle v,w \rangle = 1$.
	Reducing this identity modulo $3$,
	and using the fact that the off-diagonal entries of $g$ are multiplies of $\sqrt{-3}$, we find that $\bar g_{1,1} g_{3,3} \equiv 1 \bmod 3$.
	If $g$ is in $\Upsilon$, then this implies
	$g_{3,3}\equiv 1 \bmod 3$.
	Using the fact that $\det g = 1$, we
	obtain the other congruence $g_{2,2} \equiv 1 \bmod 3$.
	
	We'll next show that the subset
	$\Upsilon$ defined in the lemma is a subgroup.
	Suppose $g$ and $h$ are elements of
	$\Upsilon$ then
	\[
		(gh)_{1,1}
		=
		g_{1,1} h_{1,1} + g_{1,2}h_{2,1} + g_{1,3}h_{3,1}.
	\]
	We have congruences $g_{1,1} \equiv h_{1,1} \equiv 1 \bmod 3$. Furthermore 	$g_{1,2}$, $g_{1,3}$, $h_{2,1}$ and $h_{3,1}$ are all multiples of $\sqrt{-3}$.
	Therefore $(gh)_{1,1}\equiv 1 \bmod 3$, and so $gh \in \Upsilon$.
	Since $g^{-1} = J \bar g^t J$,
	we have $(g^{-1})_{1,1} = \bar g_{3,3}$.
	The congruence $g_{3,3} \equiv 1 \bmod 3$
	implies $g^{-1}\in \Upsilon$.
	
	It is clear that $\Upsilon$ and $Z$ are subgroups of $\Gamma(\sqrt{-3})$ which commute with each other and have trivial intersection,
	so it only remains to show that $\Gamma(\sqrt{-3})=\Upsilon \cdot Z$.
	To see why this is the case,
	suppose $g\in \Gamma(\sqrt{-3})$.
	We have $g_{1,1} \equiv 1 \bmod \sqrt{-3}$.
	Hence there is a cube root of unity $\zeta^r$, such that $g_{1,1} \equiv \zeta^r\bmod 3$.
	It follows that $\zeta^{-r} g \in \Upsilon$,
	so $g \in \Upsilon \cdot Z$.
\end{proof}

We shall find a presentation for $\Upsilon$.
We start by finding a small generating set.
It will be useful to have the following notation for
upper-triangular elements of $\Upsilon$:
\[
	n(z,x)
	=
	\begin{pmatrix}
		1 & \sqrt{-3} \cdot z & \frac{-3\Norm(z) + x\sqrt{-3}}{2}\\[2mm]
		0 & 1 & \sqrt{-3}\cdot\bar z\\[2mm]
		0 & 0 & 1
	\end{pmatrix},
	\qquad
	z\in \Z[\zeta],\quad x\in\Z,\quad x\equiv \Norm(z) \bmod 2.
\]
Here and later we write $\Norm$ and $\Tr$ for the norm and trace maps from
$\Q(\zeta)$ to $\Q$.

\begin{lemma}
	\label{lemma:Euclidean algorithm}
	The group $\Upsilon$ is generated by the elements $n(z,x)$ and their transposes $n(z,x)^t$, where $z\in \Z[\zeta]$ and $x\in\Z$ satisfy $x \equiv \Norm (z) \bmod 2$.
\end{lemma}

\begin{proof}	
	Let $g$ be an element of $\Upsilon$, and let
	\[
		g = \begin{pmatrix}
			a & * & * \\
			b & * & * \\
			c & * & *
		\end{pmatrix},
		\qquad
		a \equiv 1 \bmod 3,\quad
		b \equiv c \equiv 0 \bmod \sqrt{-3}.
	\]
	Define $H(g) = \Norm(a)+\Norm(c)$.
	Since $g \in \SU(2,1)$ we have $\Tr(a \bar c) + \Norm (b) = 0$.
	In particular, $a$ and $c$ cannot both be zero,
	so $H(g)$ is a positive integer.
	We shall prove by induction on $H(g)$,
	that $g$ may be expressed as a product of elements of the form $n(z,x)$ and $n(z,x)^t$.
	
	Assume first that $H(g)=1$.
	Since $\Norm(c)$ is a multiple of $3$, we must have $\Norm(a)=1$ and $c=0$. This implies $b=0$.
	The congruence $a \equiv 1 \bmod 3$ implies $a=1$. Using the fact that $g \in \SU(2,1)$ we deduce that $g= n(z,x)$ for suitable $z$ and $x$.
	
	Assume now that $H(g) > 1$.
	Our congruence conditions on $a$ and $c$ imply that $\Norm(a) \ne \Norm(c)$.
	For the inductive step, we must prove the following:
	\begin{enumerate}
		\item
		If $\Norm(a) > \Norm(c)$
		then there exists an element $n(z,x)$ such that $H(n(z,x) g) < H(g)$.
		\item
		If $\Norm(a) < \Norm(c)$
		then there exists an element $n(z,x)^t$ such that $H(n(z,x)^t g) < H(g)$.
	\end{enumerate}
	We shall prove statement 1 in detail; the proof of statement 2 is similar.
	
	Assume that $\Norm(a) > \Norm(c)$.
	One may check that the closed hexagon in $\C$	with vertices at the roots of unity $\pm 1,\pm \zeta,\pm\zeta^2$ is a fundamental domain for the lattice $\sqrt{-3}\cdot \Z[\zeta]$.
	If we choose $z\in \Z[\zeta]$ so that $\sqrt{-3}\bar z$ is
	as near as possible to $-\frac{b}{c}$,
	then $\frac{b}{c}+\sqrt{-3}\bar z$ is in this closed hexagon.
	In particular $\Norm(\frac{b}{c}+\sqrt{-3}\bar z) \le 1$.
	We have
	\[
		n(z,x_0)\cdot g
		=\begin{pmatrix}
			a' & * & * \\
			b' & * & * \\
			c & * & *
		\end{pmatrix},
		\qquad
		\text{where }
		b' = b + \sqrt{-3}\bar z c.
	\]
	Here we have chosen $x_0$ to be any integer congruent to $\Norm(z)$ modulo $2$; this choice does not change $b'$.
	Our choice of $z$ implies
	\begin{equation}
		\label{eqn:norm bound}
		\Norm\left(\frac{b'}{c}\right) \le 1.
	\end{equation}
	Next we choose $x\in \Z$ so that
	$\sqrt{3}x$ is as near as possible to
	$\Im(\frac{a'}{c})$.
	With this choice we have
	$|\Im(\frac{a'}{c} - \sqrt{-3}x)| \le \frac{\sqrt{3}}{2}$.
	Consider the matrix
	\[
		n(z,x_0-2x)\cdot g
		=\begin{pmatrix}
			a'' & * & * \\
			b' & * & * \\
			c & * & *
		\end{pmatrix},
		\qquad
		\text{where }
		a'' = a' - \sqrt{-3} x c.
	\]
	Since this matrix is in $\SU(2,1)$ we must have $a'' \bar c + b' \bar b' + c \bar a'' = 0$.
	This implies $\Re( \frac{a''}{c}) = -\frac{1}{2}\Norm(\frac{b'}{c})$, so by
	(\ref{eqn:norm bound}) we have
	$\left|\Re( \frac{a''}{c})\right| \le \frac{1}{2}$.
	Our choice of $x$ implies $|\Im(\frac{a''}{c})| \le \frac{\sqrt{3}}{2}$.
	Combining these bounds, we obtain
	\[
		\Norm\left(\frac{a''}{c}\right)
		= \Re \left(\frac{a''}{c}\right)^2 + \Im\left(\frac{a''}{c}\right)^2 \le 1.
	\]
	In particular, $\Norm(a'') \le \Norm(c) < \Norm(a)$ and therefore $H(n(z,x_0-2x) \cdot g)<H(g)$.
\end{proof}

\subsection{The weight denominators of $\Upsilon$ and $\Gamma(\sqrt{-3})$}

To calculate the weight denominators of $\Upsilon$ and $\Gamma(\sqrt{-3})$ we shall use a presentation of $\Upsilon$.

\begin{proposition}
	\label{prop:Upsilon presentation}
	The group $\Upsilon$ has a presentation with
	the following five generators:
	\[
		n_1=n(1,1),	\quad
		n_2=n(\zeta,1),\quad
		n_3=n(0,2),\quad
		n_4 = n_1^t,\quad
		n_5 = n_3^t,
	\]
	and the following thirteen relations:
	\[
	\begin{matrix}
		[n_1,n_3], \qquad
		[n_2, n_3], \qquad
		[n_4, n_5],\qquad
		( n_3  n_5 )^3,\\[2mm]
	    n_3 n_2 n_1 n_2^{-1} n_3 n_1^{-1} n_3,\qquad\qquad
    	(n_1^{-1} n_3 n_4^{-1})^3,\\[2mm]
    	n_5^{-1} n_2 n_5 n_4^{-1} n_1^{-1} n_2^{-1} n_3 n_4 n_3^{-1} n_1,\\[2mm]
	    n_4^{-1} n_1^{-1} n_3 n_5 n_2 n_1 n_5^{-1} n_4 n_2^{-1} n_3^{-1},\\[2mm]
	    n_5^{-1} n_4 n_1 n_5 n_3^{-1} n_1^{-1} n_2^{-1} n_4^{-1} n_3^{-2} n_2,\\[2mm]
	    n_5^{-1}n_2 n_1 n_5^{-1} (n_4 n_1)^2 n_5^{-1} n_4 n_2^{-1},\\[2mm]
	    n_3 n_5 n_1 n_4 n_5^{-1} n_2^{-1} n_4 n_3^{-1} n_1 n_4 n_2 n_1,\\[2mm]
		n_3^{-1} n_1 n_4 n_2 n_3 n_1 n_5^{-1} n_1^{-1} n_4^{-1} n_5 n_1^{-1} n_2^{-1},\\[2mm]
		n_4^{-1} n_3^{-1} n_5 n_3 n_1^{-1} n_4^{-1} n_2
		n_1 n_3^{-1} n_4 n_1 n_5^{-1} n_4 n_2^{-1} n_1^{-1} n_3.
    \end{matrix}
	\]
\end{proposition}

\begin{proof}
	We begin by showing that the matrices $n_1,\ldots,n_5$ generate $\Upsilon$.
	It is easy to see that any matrix $n(z,x)$ in
	$\Upsilon$ may be expressed as a product of the elements $n_1,n_2,n_3$ raised to appropriate powers.
	Similarly, any matrix $n(z,x)^t$ may be expressed in terms of $n_4,n_5, n_2^t$.
	In view of \autoref{lemma:Euclidean algorithm}, this shows that $\Upsilon$ is generated by $n_1,n_2,n_3,n_4,n_5,n_2^t$.
	The generator $n_2^t$ may be eliminated using the relation
	$n_2^t =n_3^{-1} n_1 n_4 n_1 n_3^{-2} n_2$.

	Finding the relations in $\Upsilon$ is much harder, and we have used a computer for this (see \cite{hill code} for the code).
	We'll explain briefly how the calculation was done.
	Consider the following groups:
	\begin{align*}
		\U(2,1)(\Z[\zeta])
		&=
		\{ g \in \GL_3(\Z[\zeta]) :
		\bar g^t J g = J \}, &
		\PU(2,1)(\Z[\zeta])
		&=
		\U(2,1)(\Z[\zeta]) / Z_6,
	\end{align*}
	where $Z_6$ is the centre of $\U(2,1)(\Z[\zeta])$, which is a cyclic group of order $6$.
	In \cite{falbel parker}, Falbel and Parker have obtained a presentation for the group $\PU(2,1)(\Z[\zeta])$ with three generators and five relations.
	By the method described in \autoref{sec:covering groups} we may use their result to obtain a presentation of $\U(2,1)(\Z[\zeta])$.
	The group $\Upsilon$ is a subgroup of $\U(2,1)(\Z[\zeta])$ of finite index, so we may use the Reidemeister--Schreier algorithm described
	in \autoref{sec:Reidemeister Schreier}
	to obtain a presentation of $\Upsilon$.	
	The resulting presentation of $\Upsilon$ is rather large (52 generators and 223 relations).
	The next step is to use the method described in \autoref{sec:rewriting} to rewrite the presentation in terms of the five generators $n_1,n_2,n_3,n_4,n_5$ found above.
	Finally, we simplify the presentation
	using GAP; this reduces the number of relations to the thirteen given above.
\end{proof}

\begin{theorem}
	\label{thm:wt denominator}
	The groups $\Upsilon$ and $\Gamma(\sqrt{-3})$ have weight denominator $1$.
\end{theorem}

\begin{proof}
These two groups are the same modulo their centre, so by \autoref{prop:wt denominator facts} it is sufficient to prove the result for $\Upsilon$.

We must show that $(I_3,1)$ is in the commutator subgroup $[\tilde\Upsilon,\tilde\Upsilon]$.
We shall write $r_1,\ldots, r_{13}$
 for the thirteen relations in $\Upsilon$ found in the previous result.
We shall regard relations as elements of the free group on the generators $n_1,\ldots, n_5$.

The relation $R = r_4^{-1} r_9^{-1} r_{10}^{-1} r_{11}$ is given by:
\begin{align*}
	R = \; &
	(n_5^{-1} n_3^{-1})^3 n_2^{-1} n_3^2 n_4 n_2 n_1 n_3
	 n_5^{-1}
	 n_1^{-1} n_4^{-1} n_5 n_2 n_4^{-1} n_5
	(n_1^{-1} n_4^{-1})^2\\
	& \quad n_5 n_1^{-1} n_2^{-1} n_5
	 n_3 n_5 n_1 n_4 n_5^{-1} n_2^{-1} n_4 n_3^{-1} n_1
	 n_4 n_2 n_1 .
\end{align*}
Notice that the relation $R$ has an interesting property: the number of times that each generator $n_i$ occurs in $R$ is equal to the number of times that that $n_i^{-1}$ occurs.
Equivalently $R$ is in the commutator subgroup of the free group on $n_1,\ldots, n_5$.

For each generator $n_i$, we may choose a lift $\hat n_i$ to $\tilde\Upsilon$.
Let $\hat R$ be the element of $\tilde\Upsilon$ which we obtain by replacing each generator $n_i$ in the word $R$ by the lift $\hat n_i$.
In fact $\hat R$ does not depend on the choices of lift $\hat n_i$ since the total degree of $n_i$ in $R$ is zero.
A short calculation using a computer shows that $\hat R= (I_3,-1)$.
Furthermore, $\hat R$ is in the commutator subgroup $[\tilde \Upsilon,\tilde\Upsilon]$, since every occurrence of a generator $\hat n_i$ is balanced out be an occurrence of $\hat n_i^{-1}$.
\end{proof}

Although \autoref{thm:wt denominator} was found using a computer, it would be possible to verify the result by hand (given a few days).
One needs only to calculate the element $\hat R$ in $\tilde{\Upsilon}$, and note that this element is in the commutator subgroup.
It would however be difficult to find the relation $R$ without the presentation of
\autoref{prop:Upsilon presentation}.

\subsection{Subgroups of $\Upsilon$}

There is code (see \cite{hill code}) which calculates the number $\denom (\Gamma)$ for certain arithmetic subgroups $\Gamma$ of $\Upsilon$.
The method of calculation is as follows.
\begin{enumerate}
	\item
	We already have a presentation of $\Upsilon$ (\autoref{prop:Upsilon presentation}).
	Using the Reidemeister--Schreier algorithm, we  find a presentation of $\Gamma$:
	\[
		\Gamma = \langle \gamma_1 , \ldots, \gamma_r | \delta_1 =\cdots = \delta_s = 1 \rangle.
	\]
	\item
	For each generator $\gamma_i$ we let $\hat\gamma_i = (\gamma_i,0) \in \tilde\Gamma$.
	For each relation $\delta_j$, we let $\hat \delta_j$ be the element of $\tilde\Gamma$ obtained by replacing each generator $\gamma_i$ in $\delta_j$ by its lift $\hat\gamma_i$.
	In $\tilde\Gamma$ we have relations
	\begin{equation}
		\label{eqn:delta hat}
		\hat \delta_j \cdot (I_3,n_j) = 1,
		\qquad
		n_j \in \Z.
	\end{equation}
	The integers $n_j$ are calculated
	by multiplying in $\tilde{\Gamma}$ using
	(\ref{eqn:sigma multiplication}).
	\item
	To each of the relations (\ref{eqn:delta hat}) we form a row vector
	in $\Z^{r+1}$, where the first $r$ entries
	are the multiplicities of the generators $\hat\gamma_i$ and the last entry is $n_j$.
If we let $M$ be the matrix formed of these rows, then we have a presentation
	\[
		\tilde\Gamma^{\ab}
		\cong \Z^{r+1} / (\Z^s \cdot M).
	\]
	\item
	The subgroup $\Z^s \cdot M$ is unchanged by integral row operations on $M$. We may transform $M$ by a sequence of such operations to a matrix $M'$ in Hermite normal form.
	The last non-zero row of $M'$ must have the form
	\[
		\begin{pmatrix}
			0 & \cdots & 0 & n
		\end{pmatrix},
		\qquad
		n > 0.
	\]
	Indeed if this were not the case, then the generator $(I_3,1)$ would have infinite order in $\tilde\Gamma^{\ab}$, contradicting \autoref{thm:finiteness}.
	\item
	The positive integer $n$ is the order of $(I_3,1)$ in $\tilde{\Gamma}^{\ab}$, which is by definition the weight denominator of $\Gamma$.
\end{enumerate}
The speed of the computation depends on the index of the subgroup $\Gamma$ in $\Upsilon$. The current version of the code is able to handle subgroups of index up to around 500 on a home computer.
The most time consuming step is currently the
row reduction in step 4. There is an opportunity (using \autoref{prop:wt denominator facts}) to speed up this step by row-reducing modulo the prime powers which divide the index of $\Gamma$ in $\Upsilon$, rather than row-reducing over $\Z$.
There is also scope for simplifying the presentation of $\Gamma$ at the end of step 1, using a smaller set of generators (as we did for $\Upsilon$ in \autoref{prop:Upsilon presentation}).
At the moment, this is not implemented.

\subsection{Congruence subgroups with weight denominator $3$}
\label{sec:index 3 subgroups}

The principal congruence subgroup $\Gamma(3)$ is a normal subgroup of $\Upsilon$ with index $81$.
One easily checks using the code at \cite{hill code} that $\Gamma(3)$ has weight denominator $3$.
This means that $\Gamma(3)$ has a multiplier system
with weight $\frac{1}{3}$.
If any reader would like to run this code and experiment further, then there are some instructions online. The sage code for this particular calculation is
\begin{verbatim}
        sage: load("SU21.sage")
        sage: G = Gamma(3)
        sage: print(G.weight_denominator())
\end{verbatim}
However, one does not need to go as far as $\Gamma(3)$ to find such a multiplier system.
We shall discuss the intermediate groups between $\Upsilon$ and $\Gamma(3)$

\begin{lemma}
	There is an isomorphism $F : \Upsilon / \Gamma(3) \to \F_3^4$.
	given by
	\[
		F\begin{pmatrix}
			g_{1,1} & g_{1,2} & g_{1,3} \\
			g_{2,1} & g_{2,2} & g_{2,3} \\
			g_{3,1} & g_{3,2} & g_{3,3}
		\end{pmatrix}
		=
		\begin{pmatrix}
			\frac{g_{1,2}}{\sqrt{-3}} & \frac{g_{1,3}}{\sqrt{-3}} &
			\frac{g_{2,1}}{\sqrt{-3}} & \frac{g_{3,1}}{\sqrt{-3}}
		\end{pmatrix}.
	\]
	(Here we are regarding the numbers $\frac{g_{i,j}}{\sqrt{-3}} \in \Z[\zeta]$ as elements of $\F_3$ by reducing modulo $\sqrt{-3}$.)
\end{lemma}

\begin{proof}
	We'll first show that the formula for $F$ gives a homomorphism from $\Upsilon$ to $\F_3^4$.
	Suppose $g,h\in \Upsilon$.
	We therefore have $g=I+\sqrt{-3}X$ and $h = I+\sqrt{-3}Y$ for
	suitable matrices $X,Y \in M_3(\Z[\zeta])$.
	This implies $gh \equiv I + \sqrt{-3}(X+Y) \mod 3$.
	Since $F(gh)$ depends only on $gh$ modulo $3$, it follows easily that $F(gh)=F(g)+F(h)$.
	
	Note that we have
	\[
		F(n(z,x)) = \begin{pmatrix}
			z & \frac{x}{2} & 0 & 0
		\end{pmatrix},
		\qquad
		F(n(z,x)^t) = \begin{pmatrix}
			0&0&z & \frac{x}{2}
		\end{pmatrix}.
	\]
	This shows that $F$ is surjective.
	
	If $g \in \Gamma(3)$ then clearly $g \in \ker(F)$.
	Conversely, suppose $g \in \Upsilon$ satifies $F(g)=0$.
	Since $g\in \Upsilon$, the diagonal entries of $g$ must be congruent to $1$ modulo $3$
	(by \autoref{lem:upsilon direct sum}).
	Since $F(g)=0$, we have the congruence
	\[
		g \equiv \begin{pmatrix}
			1 & 0 & 0 \\ 0 & 1 & \sqrt{-3} x\\ 0 & \sqrt{-3} y & 1
		\end{pmatrix}
		\mod 3.
	\]
	Substituting this congruence into the equation $\bar g^t J g = J$, we obtain $\sqrt{-3} x \equiv \sqrt{-3} y \equiv 0 \bmod 3$,
	which implies $g\in \Gamma(3)$.
\end{proof}

In view of the lemma, the intermediate group between $\Upsilon$ and $\Gamma(3)$ correspond to subgroups of $\F_3^4$.
In particular, for $a,b,c,d\in \F_3$ we define a subgroup of $\Upsilon$ by
\[
	\Upsilon_{\mathrm{index}\ 3}(a, b,c,d)
	=
	\{ g \in \Upsilon :
	a\cdot g_{1,2}+b\cdot g_{1,3}+c\cdot g_{2,1}+d\cdot g_{3,1}
	\equiv 0 \bmod 3\}.
\]
If $v$ is a non-zero vector in $\F_3^4$ then $\Upsilon_{\mathrm{index}\ 3}(v)$ has index $3$ in $\Upsilon$.
Every subgroup of index $3$ in $\Upsilon$ containing $\Gamma(3)$ has this form.
Two of these subgroups $\Upsilon_{\mathrm{index}\ 3}(v)$ and $\Upsilon_{\mathrm{index}\ 3}(v')$ are equal if and only if $v \equiv \pm v'\bmod 3$.
We therefore have 40 subgroups of index $3$ containing $\Gamma(3)$.
By \autoref{prop:wt denominator facts},
each of these subgroups has weight denominator either $1$ or $3$.
We have calculated all of these weight denominators.
For example, to calculate the weight denominator
of $\Upsilon_{\mathrm{index}\ 3}(0,0,1,0)$,
the command is
\begin{verbatim}
        sage: load("SU21.sage")
        sage: G = Index3congruence(0,0,1,0)
        sage: print(G.weight_denominator())
\end{verbatim}
The results of our calculations are as follows.

\begin{theorem}
	\label{thm:list}
	The following thirteen groups have weight denominator $3$:
	\[
		\begin{array}{c}
		\Upsilon_{\mathrm{index}\ 3}(0,0,1,0),\quad
		\Upsilon_{\mathrm{index}\ 3}(0,0,1,1),\quad
		\Upsilon_{\mathrm{index}\ 3}(0,0,1,2),\quad
		\Upsilon_{\mathrm{index}\ 3}(0,1,1,0),\\[2mm]
		\Upsilon_{\mathrm{index}\ 3}(0,1,2,0),\quad
		\Upsilon_{\mathrm{index}\ 3}(1,0,0,0),\quad
		\Upsilon_{\mathrm{index}\ 3}(1,0,0,1),\\[2mm]
		\Upsilon_{\mathrm{index}\ 3}(1,0,0,2),\quad
		\Upsilon_{\mathrm{index}\ 3}(1,0,2,0),\quad
		\Upsilon_{\mathrm{index}\ 3}(1,1,0,0),\\[2mm]
		\Upsilon_{\mathrm{index}\ 3}(1,1,2,2),\quad
		\Upsilon_{\mathrm{index}\ 3}(1,2,0,0),\quad
		\Upsilon_{\mathrm{index}\ 3}(1,2,2,1).
		\end{array}
	\]
	The other twenty seven groups of the form
	$\Upsilon_{\mathrm{index}\ 3}(v)$
	all have weight denominator $1$.
\end{theorem}

Note that for each of the groups $\Gamma$ listed above, the group $\Gamma\cdot Z$ is a
subgroup of index $3$ in $\Gamma(\sqrt{-3})$,
and also has weight denominator $3$ by
\autoref{prop:wt denominator facts}.

\section{Modular forms}
\label{section: modular forms}

Multiplier systems arise when one wants to consider modular forms of non integral weight.
It is a natural problem to realize a given multiplier system by a modular form. In some cases the
whole rings of modular forms of integral weight have been determined and we can use these results
to construct forms of non integral weight. Recall that a modular form on an arithmetic
subgroup $\Gamma\subset \U(n,1)$ and with respect to a multiplier system
$\ell(g,\tau)$
is a holomorphic function $f:\cH\to\C$ with the property
$$f(g\tau)=\ell(g,\tau)f(\tau),\quad g\in\Gamma,$$
where in the case $n=1$ the usual regularity condition at the cusps has to be added.
\smallskip
Notice that we work in this section with the group $\U(n,1)$ instead of $\SU(n,1)$.
This is due to the fact that in the theory of modular forms reflections play an important role.
Of course
Definition 1 works literally in the case
$\U(n,1)$.
A general result using Poincar\'e series or compactification theory states that
for every multiplier system of weight $\frac{a}{b}$ there exists a natural number $k\gg0$ such that
the multiplier system $\ell(g,\tau)j(g,\tau)^k$ of weight $\frac{a}{b}+k$ admits a non zero modular form.
We will see that 12 of the above 13 multiplier systems
admit a modular form of weight $\frac{1}{3}$.

\subsection{Rings of modular forms}

We consider the $(n+1)\times (n+1)$ matrix
$$S=\begin{pmatrix}0&-1&&\\ -1&0&&\\&&1&\\&&&\ddots\\&&&&1\end{pmatrix}.$$
It defines a Hermitian form of signature $(n,1)$,
$$\langle z,w\rangle =\bar z^t S w,\quad z,w\in\C^{n+1}\ \hbox{(columns)}.$$
We denote by $G_n$ the group of all  $(n+1)\times(n+1)$-matrices $g$ with coefficients in the
ring of Eisenstein integers
with the property
$\bar g^t S g=S$. Notice that we do not assume that $\det(g)=1$.
\smallskip
In the case $n=2$ it is related to the form $J$ from section~\ref{sec:intro}
by
$$J=\bar V^t SV,\quad V=\begin{pmatrix}-1&0&0\\0&0&1\\0&1&0\\ \end{pmatrix}.$$
So the group $\Gamma(1)$ of section~\ref{sec:results} is embedded into $G_2$,
$$\Gamma(1)\longrightarrow G_2,\quad g\longmapsto h=VgV^{-1}.$$
For each non zero Eisenstein integer $\beta$ we can consider the principal congruence
subgroups $\Gamma(\beta)$ and $G_n(\beta)$. The group $\Gamma(\beta)$ can be identified with the
subgroup of $G_2(\beta)$ of matrices with determinant one. Rings of modular forms of integral weight
have been determined in the literature for $G_4(\sqrt{-3})$, \cite{freitag},
\cite{allcock freitag}, for $G_3(\sqrt{-3})$, \cite{freitag salvati},
\cite{kondo} and for the group $G_3(3)$ \cite{freitag salvati}.

The ring of modular form for $G_3(3)$
is rather complicated. We can use its structure to derive the rings for $G_2(2)$
(and also for $G_1(3)$). We formulate now these results and give some hints how to
get them.

First we recall the definition of the ring of modular forms. Here we use the unitary group for
the Hermitian form $S$ and denote it from now on by $\U(1,n)$. So $G_n(\ell)$ is a subgroup of
$\U(1,n)$. We have to replace the homogeneous space $\cH$ by
the space $\cH_n$ of all columns $z\in \C^n$ with the property
$${\begin{pmatrix}1\\\bar z\end{pmatrix}}^tS{\begin{pmatrix}1\\\ z\end{pmatrix}}<0.$$
This means
$$\cH_n=\{z\in\C^n;\quad \Re z_1>\vert z_2\vert^2+\cdots+\vert z_n\vert^2\}.$$
The action of $\U(n,1)$ is given by
\[
	g*z
	=\frac{1}{A+Bz} \cdot (C+Dz).
\]
Here
$g$ is a block matrix
$g=\begin{pmatrix} A&B\\C&D \end{pmatrix}$,
where $A$ is a complex number, $D$ a $2\times2$-matrix etc.
We also define
\begin{equation}
	\label{eqn: j definition 2} 
	j(g,z)
	= A + Bz.
\end{equation}
The definition of a multiplier system is the same as in Definition 1 and the notion of a modular
form has to be defined as in the beginning of section~\ref{section: modular forms}.
The ring of modular forms $A(G_n(l))$ is the direct sum of the spaces of all modular forms with respect
to integral powers of $j(g,z)$.
From compactification theory one knows that this is a finitely generated algebra whose associated
projective variety is a compactification (by finitely many points) of $\cH_n/G_n(\beta)$.
The embedding
$$\cH_n\longrightarrow \cH_{n+1},\quad z\longmapsto (z,0),$$ induces a homomorphism
$$A(G_{n+1}(\beta))\longrightarrow A(G_n(\beta)).$$
From compactification theory follows that $A(G_n(\beta))$ is the normalization of the image.

\subsection{Examples of rings of modular forms}

We apply the discussion above to determine the ring of modular forms for $G_2(3)$.
We want to use the known ring for $G_3(3)$.
This ring is rather complicated. There are 15 basic modular forms $B_1,\dots,B_{15}$ of weight 1.
They have been constructed in \cite{freitag salvati} as Borcherds products.
The algebra $\C[B_1,\dots,B_{15}]$ agrees with $A(G_3(3))$ in weight $\ge 7$.
The zeros of the $B_i$ are located at Heegner divisors which can
be considered as $G_3(\sqrt{-3})$-translates of $\cH_2$ which can  be identified with the subspace of
$\cH_3$ defined by $z_3=0$.
A finite system of generators of the ideal of relations has been determined \cite{freitag salvati}.
Using these results  it is easy to prove that the
restrictions of
$$B_1,B_2,B_3,B_4+B_5,B_6+B_7,B_8-B_9,B_{10}-B_{11},B_{12}+B_{13},B_{14}+B_{15}$$
to $\cH_2$ are zero.
The quotient of $\C[B_1,\dots,B_{15}]$ by the ideal that is generated by these relations turns out to be
normal and of Krull dimension $3$. Hence this ideal is the kernel of the
homomorphism $\C[B_1,\dots,B_{15}]\to A(G_2(3))$ and this homomorphism must be surjective. In this way we can determine
$A(G_2(3))$.
We denote the restrictions of
$$B_5,\ B_{11},\  B_{15},\ B_7,\ B_{13},\ B_9$$
to $\cH_2$ in the same order by
$$X_1,\ X_2,\ X_3,\ X_4,\ X_5,\ X_6.$$

\begin{proposition}
The algebra $A(G_2(3))$ is generated by $X_1,\dots, X_6$. Defining relations
are
$$X_1^3+X_2^3-X_3^3=0,\ X_4^3-X_2^3-X_5^3=0,\ X_6^3-X_4^3+X_3^3=0.$$
The group $G_2$ acts on the forms $X_i$ by permutations combined with muliplication with 6th roots
of unity.
\end{proposition}

As we mentioned these results are consequences of the paper [5] where the case $G_3$ has been treated.
For the action of $G_3$ on the $B_i$ we refer to
Lemma 8.3 in this paper. The action of $G_2$ on the $X_i$ is a consequence.

\subsection{Modular forms of weight $\frac{1}{3}$}

Let $f\in A(G_2(3))$ be a modular form of weight $r$. Since the structure of $A(G_2(3))$ is
known and simple enough, one can determine the primary decomposition
of the principal ideal $(f)$,
$$(f)=fA(G_2(3))=\mathfrak{q}_1\cap\dots \cap \mathfrak{q}_m.$$
If one has luck, the multiplicities of the $\mathfrak{q}_i$ are multiples of 3. Then the multiplicities
of the zero components of $f$ in $\cH_2$ are also multiples of $3$. So wa can take a
holomorphic cube root of $f$ to produce a form of weight $\frac{r}{3}$.
Hence our construction of forms of weight $1/3$ rests on the knowledge of the structure of the ring
of modular forms of {\it integral\/} weight.
We give an example.

We consider the form $X_1+X_2$ of weight one. It has been found by trial and error.
Its primary decomposition can be computed by means of MAGMA. It turns out that $(X_1+X_2)$ is the intersection
of three primary ideals

\begin{align}
\mathfrak{q}_1&=(X_1 + X_2,X_4+(\zeta+1)X_6,X_3^3,X_2^3 + X_5^3 - X_6^3),\\
\mathfrak{q}_2&=(X_1 + X_2,X_4-\zeta X_6,X_3^3,X_2^3 + X_5^3 - X_6^3),\\
\mathfrak{q}_3&=(X_1 + X_2,X_4-X_6,X_3^3,X_2^3 + X_5^3 - X_6^3).
\end{align}

The associated prime ideals (which are the radicals) $\mathfrak{p}_1,\mathfrak{p}_2,\mathfrak{p}_3$
are obtained if one replaces in each case $X_3^3$ by $X_3$.
For example
$$\mathfrak{p}_3=(X_1 + X_2,X_4-X_6,X_3,X_2^3 + X_5^3 - X_6^3).$$
This defines the elliptic curve
$$X_2^3 + X_5^3 - X_6^3=0$$
in $P^2\C$. Its $j$-invariant is $0$. All three curves are isomorphic.

\begin{lemma}
The form $X_1+X_2$ on $\cH_2$ with respect to the group $G_2(3)$
vanishes along three elliptic curves with multiplicity three and has no other
zero. Hence it is the third power of a modular form of weight $\frac{1}{3}$.
\end{lemma}

In \cite{freitag salvati} the action of $G_3$ on $A_3(G_3)$ has been determined.
Using this we can determine the invariance group of $X_1+X_2$ to verify the following result.

\begin{lemma}
The invariance froup of the form $X_1+X_2$ is
the subgroup of $G_2(\sqrt{-3})$ that is defined through the congruence
$h_{13}+h_{21}\equiv 0\mod 3$.
This group is an extension of index 3 of its intersection
with $\SU(2,1)$.
This intersection corresponds to the group $\Upsilon_{\mathrm{index}\ 3}(0,1,1,0)$
(subsection~\ref{sec:index 3 subgroups}).
\end{lemma}

The form $X_1+X_2$ is a special example of a form of weight 1 which admits a holomorphic third
root. To get more such forms, we tried also the forms $X_i\pm X_j$. It turned out that for each
pair $(i,j)$, $i\ne j$, there is one of the two signs such that the form has third root of unity.

The following table contains forms  $X_i\pm X_j$ which have the same property.
They belong to a subgroup defined by a congruence $L(h)\equiv0\mod 3$. The linear form
$L$ is in the first column. The second column contains the form and the third column contains
the vector $v$ such that intersection with $\SU(2,1)$ corresponds to
$\Upsilon_{\mathrm{index}\ 3}(v)$.
\bigskip\noindent

{\bf 12 Congruence groups that admit a modular form of weight $\frac{1}{3}$}

\halign{\qquad\qquad$#$\qquad\hfill&$#$\hfill&$\qquad#$\hfill\cr
\noalign{\smallskip}
\hbox{Congruence $L$}&\hbox{third power of a}&\Upsilon_{\mathrm{index}\ 3}(v)\cr
&\hbox{form of weight}\ 1/3&v:=\cr
\noalign{\smallskip}
h_{13}& X_2-X_3&(0,0,1,0)\cr
h_{12}+h_{13}& X_3-X_4&(0,0,1,1)\cr
h_{12}-h_{13}& X_2-X_4&(0,0,1,2)\cr
h_{13}+h_{21}& X_1+X_2&(0,1,1,0)\cr
h_{13}-h_{21}& X_1-X_3&(0,1,2,0)\cr
h_{31}& X_5-X_6&(1,0,0,0)\cr
h_{12}+h_{31}& X_1-X_5&(1,0,0,1)\cr
h_{12}-h_{31}& X_1+X_6&(1,0,0,2)\cr
h_{21}+h_{31}& X_4-X_5&(1,1,0,0)\cr
h_{12}+h_{13}-h_{21}-h_{31}& X_4-X_6&(1,1,2,2)\cr
h_{21}-h_{31}& X_3+X_6&(1,2,0,0)\cr
-h_{12}+h_{13}+h_{21}-h_{31}&X_2+X_5 &(1,2,2,1)\cr
}
\medskip
These forms are not uniquely determined. In fact, one can replace each form $X_i\pm X_j$
in the table by $X_i\pm\zeta^\nu X_j$,  $0\le\nu\le2$ and in this way we get forms with respect to the same group
that admit also third roots of unity.

This means that we constructed 36 forms of weight $\frac{1}{3}$.
So we have proved that 12 of the above 13 groups in subsection~\ref{sec:index 3 subgroups}
admit three modular forms of weight $\frac{1}{3}$.

One can show that this is false for the remaining group
which is given through $h_{13}\equiv h_{31}\mod 3$.
Its intersection with $\SU(2,1)$ corresponds to $v=(1,0,2,0)$. In this case one can prove that
no form  of weight one which belongs to this congruence group and which admits a holomorphic
third root can exist. By general arguments there must exist in $A(G_2(3))$
a modular of  weight $3r+1$,
$r\in\Z$ suitable, which admits a holomorphic third root of unity.

So far we do not know an example.

\end{document}